\documentclass{elsarticle}

\usepackage{amsmath,amssymb,amsfonts,graphicx,hyperref}
\usepackage[margin=1in]{geometry}

\usepackage{amstext}
\usepackage{stmaryrd}
\usepackage{algpseudocode}

\usepackage[normalem]{ulem}
\usepackage{graphicx}

\usepackage{array}
\usepackage{verbatim}
\usepackage{booktabs}
\usepackage{calc}
\usepackage{textcomp}
\usepackage{multirow}

\usepackage[all]{xy}

\usepackage{cancel}
\usepackage{mathtools}
\usepackage{placeins}
\usepackage{xspace}

\usepackage{tikz,tikzscale}
\usetikzlibrary{calc}
\usetikzlibrary{shadings}
\usepackage[mode=multiuser,status=draft]{fixme} % Provide note, warning, error, fatal.
\fxsetup{envlayout=color}
\fxsetup{innerlayout=inline}
\fxsetup{targetlayout=colorcb}
\FXProvidesTargetLayout{color}

\fxusetheme{color}
\FXRegisterAuthor{aa}{envmr}{AA}
\FXRegisterAuthor{mw}{envmw}{MW}

\usepackage{subcaption}
\usepackage[format=plain,indention=.5cm]{caption}
\usepackage{pgfplots}
\pgfplotsset{compat=1.18}
\usepgfplotslibrary{fillbetween}

\usepackage{colortbl} % Define a row color

\usepackage{booktabs} % Allows the use of \toprule, \midrule and \bottomrule in tables   

\usepackage{bm,upgreek}

\newcommand{\DGell}{{\mathbb{V}_\ell}}
\usepackage{graphicx} % Required for inserting images
\usepackage[linesnumbered,vlined,commentsnumbered,ruled]{algorithm2e}

\definecolor{darkgreen}{rgb}{0.0, 0.5, 0.0}

\definecolor{deepgreen}{RGB}{0,100,0}

\newcommand{\plotwidth}{0.49\columnwidth}
\newcommand{\plotheight}{0.4\columnwidth}

\pgfplotscreateplotcyclelist{paper markers}{
    {teal,    solid, mark=*},
    {orange,   solid, mark=triangle*},
    {blue,     solid, mark=square*},
    {red,      solid, mark=pentagon*},
    {violet,     solid, mark=otimes*},
    {brown,    solid, mark=halfcircle*},
    {black,   solid, mark=diamond*},
}

\pgfplotsset{cycle list name=paper markers}

\pgfplotsset{
    paperplot/.style={
            width=\plotwidth,
            height=\plotheight,
            cycle list name=paper markers,
            every axis/.append style={font=\small},
            grid=major,
            thick
        }
}

\pgfplotsset{
    longplot/.style={
            width=0.8\columnwidth,
            height=\plotheight,
            cycle list name=paper markers,
            every axis/.append style={font=\small},
            grid=major,
            thick
        }
}

\title{A Geometric Multigrid Preconditioner for Shifted Boundary Method}
\author{Micha\l{} Wichrowski\footnote{Heidelberg University, Germany, mt.wichrowsk@uw.edu.pl}, Ajay Ajith
    \footnote{Heidelberg University, Germany}}
\begin{document}

\begin{abstract}
    The Shifted Boundary Method (SBM) trades some part of the burden of body-fitted meshing for increased algebraic complexity. While the resulting linear systems retain the standard $\mathcal{O}(h^{-2})$ conditioning of second-order operators, the non-symmetry and non-local boundary coupling render them resistant to standard Algebraic Multigrid (AMG) and simple smoothers for high-order discretizations. We present a geometric multigrid preconditioner that effectively tames these systems. At its core lies the \emph{Full-Residual Shy Patch} smoother: a subspace correction strategy that  filters out some patches while capturing the full physics of the shifted boundary. Unlike previous cell-wise approaches that falter at high polynomial degrees, our method delivers  convergence with low mesh dependence. We demonstrate performance for Continuous Galerkin approximations, maintaining low and stable iteration counts up to polynomial degree $p=3$ in 3D, proving that SBM can be both geometrically flexible and algebraically efficient.
\end{abstract}
\begin{keyword}
    immersed boundary methods \sep finite element methods \sep Shifted Boundary Method \sep
    multigrid \sep  patch smoothers
\end{keyword}
\maketitle
\section{Introduction}

% Motivation: PDEs on complex geometries, unfitted methods, and SBM basics leading to conditioning issues

% Motivation: PDEs on complex geometries, unfitted methods, and SBM basics leading to conditioning issues
Capturing the intrinsic detail and complexity of the real world within the rigid structure of computational grids is a
persistent tension in scientific computing. Traditional finite element methods (FEM) require body-fitted meshes that
must conform precisely to every contour of the domain's boundary. While effective, this bespoke mesh generation is
computationally expensive and frequently becomes the bottleneck for intricate three-dimensional shapes. Unfitted finite
element methods offer a liberating alternative by employing a fixed background mesh that remains agnostic to the
physical boundary. Among these, the Shifted Boundary Method (SBM)~\cite{main2018shifted} adopts a distinct strategy:
instead of the mesh chasing the boundary, the boundary is mathematically shifted to meet the mesh. By defining the
problem on a surrogate domain and extrapolating boundary conditions, SBM avoids the complex geometric intersections of
methods like CutFEM~\cite{burman2015cutfem}. However, in a display of the "conservation of difficulty," SBM essentially
shifts part of the challenge from the mesh generator to the linear solver. The resulting linear systems are burdened
with extrapolation terms that introduce non-symmetry and potential indefiniteness, proving that a simple mesh does not
necessarily guarantee a simple matrix. While conditioning issues scale like
$\mathcal{O}(h^{-2})$~\cite{atallah2021shifted}, similar to body-fitted methods~\cite{collins2023penalty}, the
efficient numerical solution of algebraic systems emanating from high-order SBM formulations remains somewhat
unexplored.

% The gap in solvers for SBM, previous DG work limitation, and the need for patches
Despite its advantages, the development of robust and scalable solvers for SBM, particularly geometric multigrid
preconditioners~\cite{hackbusch1985multi}, has remained an open problem. In previous
work~\cite{wichrowski2025geometric}, a geometric multigrid preconditioner was developed for a Discontinuous Galerkin
(DG) SBM formulation. That approach leveraged the inherent block structure of DG to utilize element-wise smoothers.
However, while effective for low orders, that strategy faced significant difficulties at higher polynomial degrees
($p=3$), where the simple cell-wise smoother failed to resolve the error efficiently, leading to deteriorating
iteration counts. For Continuous Galerkin (CG) methods, the strong coupling between degrees of freedom makes
element-wise smoothing even less viable. Consequently, this work aims to not only address the specific challenges of
classical SBM but also to overcome the high-order robustness issues encountered in previous formulations by introducing
a more sophisticated patch-based smoothing strategy.

% Literature review: SBM details, evolution, and applications (Stokes, Interfaces, etc.)
In SBM, the \emph{surrogate domain} $\tilde{\Omega}$ is typically constructed as a union of cells from a fixed
background mesh that are deemed \emph{active} (e.g., entirely inside or significantly intersecting the true domain
$\Omega$), and its boundary $\tilde{\Gamma}$ does not conform to the domain boundary $\Gamma$. Boundary conditions are
transferred from the true to the surrogate boundary, typically via Taylor expansions or more general extension
operators~\cite{zorrilla2024shifted}, and enforced in a Nitsche-like manner~\cite{nitsche1971variationsprinzip}. The
Shifted Boundary Method has evolved from its first formulation~\cite{main2018shifted, main2018shiftedVol2}, which used
cells strictly within the considered domain, to a recent approach~\cite{yang2024optimal} that often includes
intersected cells based on a volume fraction threshold. It has been extended to high-order discretizations
\cite{atallah2022high}, various physical problems including Stokes flow~\cite{atallah2020analysis}, solid
mechanics~\cite{atallah2021shifted, atallah2024nonlinear}, and significantly, to problems with embedded
interfaces~\cite{li2020shifted, xu2024weighted}. These latter works demonstrate the capability of SBM to handle
discontinuities across internal boundaries by appropriately modifying the formulation to impose jump conditions,
extending the method's applicability to multiphysics and multi-material problems. Other innovations include
penalty-free variants~\cite{collins2023penalty} and integration with level set methods~\cite{kuzmin2022unfitted,
    xue2021new}.

% Comparison with CutFEM: Geometric flexibility vs. Quadrature efficiency
While SBM avoids the complexities of generating body-fitted meshes, the primary geometric task shifts to accurately
determining the relationship between points on the surrogate boundary $\tilde{\Gamma}$ and the true boundary $\Gamma$.
The method inherently allows for the use of arbitrarily complex geometries, and crucially, avoids the need to compute
integrals over the arbitrarily shaped integration domains that arise from cell-boundary intersections. While
determining the active mesh requires computing volume fractions of cut cells, this is a one-time geometric
preprocessing step. In contrast, methods like CutFEM require specialized quadrature for all terms in the bilinear form
on every cut cell. This makes SBM significantly more efficient in terms of computational throughput, especially in
matrix-free implementations where the overhead of cut-cell quadrature would be incurred in every operator evaluation.
The SBM utilizes closest-point projection algorithms to find for each point on $\tilde{\Gamma}$ a corresponding point
on $\Gamma$. Level sets,which represent the domain boundary as the zero level set of a function, are a common way to
represent an embedded domain, and SBM uses closest-point projection algorithms to find for each point on
$\tilde{\Gamma}$ a corresponding point on $\Gamma$. Level set methods, which represent the domain boundary as the zero
level set of a function, are often employed in unfitted methods like SBM to facilitate operations such as closest point
projection~\cite{kuzmin2022unfitted, xue2021new}. Special treatment of domains with corners was analyzed
in~\cite{atallah2021analysis}.

% Algebraic challenges: Non-symmetry and indefiniteness
However, the geometric flexibility of SBM comes at the cost of new computational challenges. The extrapolation of
boundary conditions and the resulting modifications to the variational formulation often lead to linear systems with
non-symmetric and potentially indefinite properties~\cite{wichrowski2025geometric}, particularly for higher-order
polynomial approximations. These properties make the efficient solution of the resulting system non-trivial and demand
specialized preconditioning strategies.

% Multigrid basics and failure of standard smoothers/AMG
Multigrid methods~\cite{brandt1977multi, hackbusch1985multi} are renowned for their potential to solve large systems
arising from partial differential equations, often offering optimal or near-optimal complexity. However, their
application as preconditioners for SBM remains largely unexplored. While Algebraic Multigrid (AMG) has been applied to
SBM for discretizations using continuous linear elements~\cite{atallah2020second}, its efficiency for higher-order
methods is not established. In fact, AMG generally struggles with high-order finite element discretizations even for
body-fitted meshes, and for unfitted methods like SBM, it is even less effective for $p > 1$. Furthermore, AMG may not
fully leverage the geometric information in structured background meshes. Geometric multigrid methods, in contrast,
explicitly use the hierarchy of meshes and can be very effective. A critical component is the smoother. Standard
smoothers like Jacobi or Gauss-Seidel, however, struggle with the non-standard local properties of SBM systems. They
implicitly enforce zero correction on the boundary of the local support (e.g., a vertex patch for Gauss-Seidel), which
conflicts with the non-zero boundary terms introduced by SBM.

% Solution strategy: Subspace correction and Patch smoothers
To address these issues in the context of Continuous Galerkin SBM, we adopt the perspective of subspace correction
methods~\cite{xu2001method, bramble1991analysis}. The effectiveness of a smoother depends critically on the choice of
subspaces and the local problems solved on them. Since simple one-dimensional subspaces (corresponding to individual
DoFs) are ineffective for SBM due to the boundary conflicts mentioned above, we construct local problems on larger,
overlapping subspaces defined by patches of elements centered around vertices. This concept is akin to overlapping
Schwarz smoothers~\cite{pavarino1993additive}, but tailored for the specific challenges of SBM. By solving a local
problem on a patch that includes the relevant SBM boundary terms, we can compute a more effective correction. We
introduce a \emph{shyness} criterion to ensure that patches are only formed around vertices sufficiently surrounded by
active cells, thereby preventing the creation of ill-conditioned local problems on small, isolated slivers of the
domain.

% Matrix-free context and efficiency of patch smoothers
The adoption of patch-based smoothers is further justified by recent algorithmic improvements that render them highly
efficient, particularly in matrix-free settings. While traditional overlapping Schwarz methods are often viewed as
computationally expensive, recent work on smoothers with localized residual computations~\cite{wichrowski2025smoothers}
demonstrates that multiplicative formulations allow for the fusion of residual computation with the local subspace
correction. This significantly reduces memory transfers, a key bottleneck on modern hardware. Furthermore, the local
problems on these patches can be solved with remarkable efficiency; interior patches can leverage fast tensor-product
inverses~\cite{witte2021fast, cui2025multigrid}, while more general patches can be solved using local $p$-multigrid
techniques. This approach yields $p$-robustness at a cost comparable to a global operator
evaluation~\cite{wichrowski2025multigrid, wichrowski2025local}. Indeed, recent results for CutFEM on
GPUs~\cite{cui2025multigrid} validate the effectiveness of such vertex-patch smoothers for unfitted methods, reporting
promising iteration counts and demonstrating that performing multiple smoothing sweeps specifically on boundary patches
is a viable strategy to ensure robustness. Thus, while this work utilizes matrix-based implementations to demonstrate
the robustness of the multigrid approach, the proposed smoother is designed to align with future high-performance,
matrix-free solvers that avoid expensive matrix assembly~\cite{kronbichler2019fast,wichrowski2025largeStrain,
    witte2021fast}.

% Context: CutFEM and stabilization
In the broader context of unfitted methods, CutFEM~\cite{burman2015cutfem} is a prominent alternative with theoretical
results for preconditioners already developed. It discretizes directly on the physical domain by cutting background
cells, requiring specialized quadrature. Due to inherent difficulties with small cuts, proper stabilization seems to be
an unavoidable part of CutFEM. The so-called ghost penalty~\cite{burman2010ghost, wichrowski2025matrix} solves the
issue of ill-conditioning but may require additional care to avoid
locking~\cite{badia2022linking,bergbauer2024high,burman2022design}. CutFEM has been applied to various problems,
including Stokes~\cite{burman2014fictitious}, elasticity~\cite{hansbo2017cut}, or two-phase
flows~\cite{claus2019cutfem}. Furthermore, Discontinuous Galerkin methods have also been combined with
CutFEM~\cite{gurkan2019stabilized, bergbauer2024high}. While CutFEM ensures robust conditioning via geometry-adapted
quadrature, SBM retains the efficiency comparable to standard tensor-product quadrature.

% Context: CutFEM preconditioners
Concerning preconditioning, CutFEM seems to pose challenges. Results providing optimal
preconditioners~\cite{gross2023analysis, gross2021optimal} have been developed. Although these methods are shown to be
mesh-independent, iteration counts can be high. In~\cite{bergbauer2024high} DG-CutFEM was considered, and a multigrid
preconditioner based on cell-wise Additive Schwarz smoother was used. Although the paper mostly focuses on matrix-free
implementation, the preconditioner seems promising. However, the smoothing step requires a rather high number of
matrix-vector products.

% Contribution of this paper: Full-Residual Shy Patch Smoother for CG-SBM
This paper addresses the computational challenges of solving linear systems arising from Continuous Galerkin SBM
discretizations. Our main contribution is the development of a novel subspace correction smoother, termed the
\emph{Full-Residual Shy Patch} smoother. This smoother is built on local problems defined on overlapping patches of
elements, constructed to be faithful representations of the global SBM problem. We introduce a \emph{shyness} criterion
to ensure robustness and avoid ill-conditioning on small patches. We demonstrate that this patch-based smoother,
combined with a multi-stage strategy, leads to a highly effective $h$-multigrid preconditioner. Numerical experiments
demonstrate the preconditioner's effectiveness under mesh refinement and compare our results with an algebraic
multigrid (AMG) preconditioner, highlighting the limitations of AMG for higher-order SBM discretizations. The
implementation is built on the \texttt{deal.II} finite element library~\cite{dealii2019design, dealii2024}, leveraging
its comprehensive tools for finite element methods and multigrid.

The remainder of this paper is organized as follows. In Section~\ref{sec:method}, we briefly review the Shifted
Boundary Method formulation to establish the necessary notation and context. Section~\ref{sec:mg_preconditioner}
details the core of our contribution: the construction of the Full-Residual Shy Patch smoother and the associated
multigrid hierarchy. The key concept of \emph{Shy Patches}, which ensures robustness by avoiding ill-conditioned
slivers near the boundary, is formally introduced in Section~\ref{sec:shy_patches}. We then present numerical evidence
in Section~\ref{sec:numerical_results}, demonstrating the effectiveness of the method and providing comparisons with
AMG and CutFEM. In essence, this work demonstrates that while SBM matrices can be notoriously difficult to handle at
high polynomial degrees, they can be effectively tamed by our \emph{Shy Patch} smoothers. By solving local problems on
vertex patches while carefully avoiding ill-conditioned regions, we achieve a solver that remains robust in both mesh
size $h$ and polynomial degree $p$. This robustness is clearly illustrated by the stable iteration counts presented in
Figure~\ref{fig:iteration_counts_shy3_smoothing3} and Table~\ref{tab:pmg_iterations}, while the comparative performance
against other methods is summarized in Table~\ref{tab:comparison_cutfem}. Finally, Section~\ref{sec:conclusion} offers
concluding remarks.

\section{Method formulation}
\label{sec:method}
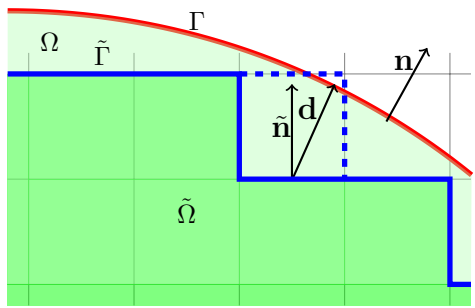
\begin{figure}[!ht]
    \centering
\begin{tikzpicture}[scale=1.4]
    % Draw the Cartesian grid
    \draw[step=1cm, gray, very thin] (0.8,-0.2) grid (5.2,2.2);

    % Highlight cells inside the domain (interior cells)
    \foreach \x/\y in {2/1, 1/0, 1/1, 2/0, 3/0, 4/0}
        {
            \fill[green, opacity=0.5] (\x,\y) rectangle (\x+1,\y+1);
        }
    \foreach \x/\y in {1/0 }{
            \fill[green, opacity=0.5] (\x,\y) rectangle (\x-0.2,\y+2);
        }

    \fill[green, opacity=0.5] (0.8,-0.2) rectangle (5.2,0);

    \draw[line width=2pt, red, domain=51:90, samples=100, variable=\t]
    plot ({0.8 + 7.0*cos(\t)}, {-4.4 + 7.0*sin(\t)});
    \node at (2.6,2.5) {$\Gamma$};

    % Fill the area below the red line with light green to indicate the domain Omega
    \fill[green!30, opacity=0.4] plot[domain=51:90, samples=100, variable=\t] ({0.8 + 7.0*cos(\t)}, {-4.4 +
            7.0*sin(\t)}) -- (0.8,0) -- (5.2,0) -- cycle;

    % Define the angle parameter
    \pgfmathsetmacro{\myangle}{66}
    % Draw a vector from the surrogate boundary to the true boundary
    \draw[<-, thick] ({0.8 + 6.9*cos(\myangle) + 0.3}, {-4.4 + 6.9*sin(\myangle)}) -- ({0.8 + 5.9*cos(\myangle)
            +
            0.3}, {-4.4 +
            5.9*sin(\myangle)});
    % Position the label relative to the vector angle
    \pgfmathsetmacro{\labelangle}{\myangle - 6} % Adjust angle slightly for label placement
    \node at ({0.8 + 6.9*cos(\labelangle) - 0.6}, {-4.4 + 6.9*sin(\labelangle)+0.1}) {\large $\mathbf{d}$};

    \draw[->, thick] ({0.8 + 6.8*cos(\myangle-5) + 0.3}, {-4.4 + 6.8*sin(\myangle-5)}) -- ({0.8 +
            7.6*cos(\myangle-5)
            +
            0.3}, {-4.4 +
            7.6*sin(\myangle-5)});
    \node at ({0.8 + 7.5*cos(\labelangle-5) - 0.54}, {-4.4 + 7.5*sin(\labelangle-5)+0.36}) {\large $\mathbf{n}$};

    % Draw the vector d vertically from the surrogate boundary to the true boundary
    \draw[<-, thick]  ({0.8 + 5.9*cos(\myangle) +
            0.3}, {-4.4 + 6.9*sin(\myangle)}) -- ({0.8 + 5.9*cos(\myangle) +
            0.3}, {-4.4 +
            5.9*sin(\myangle)});
    \node at (3.4, 1.5) {\large $\tilde{\mathbf{n}}$}; % Place label near the normal vector

    % Draw the surrogate boundary (only the top boundary of the interior cells)
    \draw[line width=2pt, blue] (0.8,2) -- (3,2) -- (3,1) -- (4,1) -- (5,1)-- (5,0)-- (5.2,0.) ;

    \draw[line width=2pt, blue, dashed] (3,2) -- (4,2) -- (4,2) -- (4,1) ;
    \node at (1.7,2.2) {$\tilde{\Gamma}$};

    \node at (2.5,0.7) {$\tilde{\Omega}$};

    \node at (1.2,2.3) {${\Omega}$};

\end{tikzpicture}
    \caption{Schematic illustrating the background mesh, interior cells (green), the surrogate boundary
        $\tilde{\Gamma}$ (thick blue line) along the upper boundary of the interior cells, and the true boundary
        $\partial\Omega$ (red). Figure taken from~\cite{wichrowski2025geometric}.}
    \label{fig:sbm_setup}
\end{figure}

We consider the Poisson problem as a model problem:
\begin{align}
    -\Delta u & = f \quad \text{in } \Omega, \label{eq:poisson}      \\
    u         & = g \quad \text{on } \Gamma, \label{eq:dirichlet_bc}
\end{align}
where $\Omega \subset \mathbb{R}^d$ ($d=2,3$) is a domain with boundary $\Gamma =\partial\Omega$ as depicted in
Figure~\ref{fig:sbm_setup}, $f$ is a given source term, and $g$ is the prescribed Dirichlet boundary condition.

To solve this problem numerically, we first formulate it in a weak sense. We seek a solution $u$ in an appropriate
function space, $H^1(\Omega)$, which consists of functions that are square-integrable and whose first derivatives are
also square-integrable. Multiplying the equation by a test function $v \in H^1(\Omega)$ and integrating over $\Omega$,
we obtain:
\[
    \int_\Omega \nabla u \cdot \nabla v \, dx - \int_{\Gamma} (\nabla u \cdot \mathbf{n}) \, v \, ds =
    \int_\Omega f \, v \, dx.
\]

We next introduce a triangulation \( \mathcal{T}_h \) consisting of quadrilateral (2D) or hexahedral (3D) elements of
size $h$, and define a finite element space $\mathbb{V}_h\subset H^1(\Omega)$ using Lagrange polynomial elements of
degree \( p \). In classical finite element methods, the mesh conforms to the boundary (unlike the background mesh
approach illustrated in Figure~\ref{fig:sbm_setup}), and test functions $v$ typically vanish on $\Gamma $ to strongly
enforce the homogeneous Dirichlet condition $u=g$ for $g=0$. When function spaces do not necessarily satisfy essential
boundary conditions strongly, boundary conditions must be enforced weakly through additional integral terms on
$\Gamma$.

Nitsche's method provides a way to weakly impose Dirichlet boundary conditions within a variational formulation without
requiring the function space to satisfy the boundary conditions. It modifies the bilinear form by adding terms on the
boundary $\Gamma$. The standard Nitsche formulation for the Poisson problem with Dirichlet boundary conditions $u=g$ on
$\Gamma$ seeks $u_h \in V_h$ such that for all $v_h \in \mathbb{V}_h$:
\begin{equation*}
    \begin{split}
        \int_\Omega \nabla u_h \cdot \nabla v_h \, dx
        \;\; - \int_{\Gamma} (\nabla u_h \cdot \mathbf{n}) v_h \, ds
        - & \alpha\int_{\Gamma} (\nabla v_h \cdot \mathbf{n}) u_h \, ds
        + \int_{\Gamma} \sigma u_h v_h \, ds =                          \\
          & = \int_\Omega f v_h \, dx
        - \alpha\int_{\Gamma} (\nabla v_h \cdot \mathbf{n}) g \, ds
        + \int_{\Gamma} \sigma_\Gamma g v_h \, ds.
    \end{split}
\end{equation*}
Here, $\sigma_\Gamma$ is a penalty parameter, typically chosen as $\sigma_\Gamma = \mathcal{O}(p^2 h^{-1})$ for mesh
size $h$, and
$\mathbf{n}$ is the outward unit normal vector to $\Gamma$.  The choice of parameter $\alpha=1$ leads to a symmetric
formulation, while $\alpha=-1$ results in a non-symmetric formulation in which the penalty term can be
skipped~\cite{baumann1999discontinuous}

\subsection{Shifted Boundary Method}
In many applications, the domain $\Omega$ may have a complex geometry, making the generation of body-fitted meshes
challenging. Non-body-fitted (unfitted) methods address this challenge by employing a background mesh $\mathcal{T}_h$
that does not conform to the boundary of $\Omega$. The domain $\Omega$ is embedded within this background mesh, and
cells of $\mathcal{T}_h$ are classified as active based on their intersection with $\Omega$. The computational domain
$\tilde{\Omega}$ is defined as the union of these active cells. In the original SBM formulation~\cite{main2018shifted},
only cells strictly contained in $\Omega$ were included, leading to the surrogate domain boundary depicted by the solid
blue line in Figure~\ref{fig:sbm_setup}. Consequently, the surrogate domain $\tilde{\Omega}$ is generally a subset of
$\Omega$, and its boundary $\tilde{\Gamma}$ does not coincide with the true boundary $\Gamma$.

Later extensions include intersected cells~\cite{yang2024optimal} with a volume fraction outside $\Omega$ less than a
threshold $\lambda\in [0,1]$. In Figure \ref{fig:sbm_setup}, this corresponds to including additional cells as
indicated by the dashed blue line. This approach reduces the distance between true and surrogate boundaries, but can
lead to ill-conditioning.

To impose boundary conditions on $\tilde{\Gamma}$, the Dirichlet condition prescribed on the true boundary $\Gamma$ is
extrapolated to the surrogate boundary. This is typically accomplished by projecting points from $\tilde{\Gamma}$ onto
the true boundary $\Gamma$ along a suitable direction (e.g., the outward normal to $\tilde{\Gamma}$) and using a Taylor
expansion to approximate the boundary values. This is accomplished via an \emph{extension operator} $\mathcal{E}$,
which maps functions defined on $\tilde{\Omega}$ to the surrogate boundary $\tilde{\Gamma}$.

We assume that the Dirichlet boundary condition $g$ is given as a restriction of a function $u^\star$ defined on the
entire domain $\Omega$ to the boundary $\Gamma$. While the choice of this function $u^\star$ (as an extension of $g$
from $\Gamma$ into $\Omega$) is not unique, we take $u^\star$ to be equal to the solution $u$ in the surrogate domain
$\tilde{\Omega}$. For each point $\tilde{\mathbf{x}} \in \tilde{\Gamma}$, let $\mathbf{x} \in \Gamma$ be its closest
point projection onto the true boundary, and let $\mathbf{d} = \mathbf{x} - \tilde{\mathbf{x}}$ be the shift vector.
The function $u$ is extended from $\Gamma$ to $\tilde{\Gamma}$ using a Taylor expansion:
\[
    \mathcal{E}u^\star(\tilde{\mathbf{x}}) = u^\star(\mathbf{x}) + \mathbf{d} \cdot \nabla u (\tilde{\mathbf{x}}) +
    \cdots
\]
where $\mathcal{E}u^\star$ denotes the extrapolated boundary condition. The function $u^\star$ is assumed to be smooth
in a neighborhood of $\tilde{\Gamma}$, which allows for the Taylor expansion to be valid. By substituting it into the
weak formulation, we obtain a variational formulation for the shifted boundary problem with the extrapolated boundary
condition on $\tilde{\Gamma}$ enforced in a Nitsche-like manner; the weak formulation seeks $u_h \in V_h$ such that for
all $v_h \in V_h$,
\begin{equation}
    \begin{split}
        \int_{\tilde{\Omega}} \nabla u_h \cdot \nabla v_h \, dx
        - \int_{\tilde{\Gamma}} (\nabla u_h \cdot \tilde{\mathbf{n}}) v_h \, ds
        - \alpha \int_{\tilde{\Gamma}} (\nabla v_h \cdot \tilde{\mathbf{n}}) \; \mathcal{E}u_h \, ds
        + \int_{\tilde{\Gamma}} \sigma \; \mathcal{E} u_h  \; v_h \, ds = \\
        = \int_{\tilde{\Omega}} f v_h \, dx
        - \alpha\int_{\tilde{\Gamma}} (\nabla v_h \cdot \tilde{\mathbf{n}}) {g} \, ds
        + \int_{\tilde{\Gamma}} \sigma g v_h \, ds.
    \end{split}
\end{equation}
where $\tilde{\mathbf{n}}$ is the outward normal to $\tilde{\Gamma}$ and $\sigma$ is a penalty parameter.

The choice of the stabilization term, particularly the parameter $\alpha$, significantly influences the spectral
properties of the resulting system matrix. As detailed in our previous work on a Discontinuous Galerkin (DG) based SBM
multigrid preconditioner~\cite{wichrowski2025geometric}, this can be observed even in a simple 1D single-cell problem.
For instance, a penalty-free formulation ($\sigma=0$, $\alpha=-1$) can lead to complex eigenvalues, especially when the
true boundary lies inside the surrogate domain.

In the context of DG methods, it was found that the choice of stabilization did not have a dramatic impact on the
overall multigrid performance, as the cell-based nature of the DG smoother effectively handled local
issues~\cite{wichrowski2025geometric}. However, in the DG formulation the cellwise smoother experienced singifficant
challenges for $p=3$. For continuous finite elements, where our smoother operate on larger, overlapping patches of
elements (i.e. vertex or edge patches), the choice of stabilization becomes far more critical. The non-standard terms
introduced by SBM are not well-contained within these patches, and an inappropriate stabilization can introduce
spectral properties that standard smoothers cannot handle effectively, a problem that is exacerbated for higher-order
elements.

The SBM weak formulation with symmetrized stabilization seeks $u_h \in V_h$ such that for all $v_h \in V_h$,
\begin{equation}
    \begin{split}
        \int_{\tilde{\Omega}} \nabla u_h \cdot \nabla v_h \, dx
        - \int_{\tilde{\Gamma}} (\nabla u_h \cdot \tilde{\mathbf{n}}) v_h \, ds
        -  \int_{\tilde{\Gamma}} (\nabla v_h \cdot \tilde{\mathbf{n}}) \; \mathcal{E}u_h \, ds
        + \int_{\tilde{\Gamma}} \sigma_\Gamma \; \mathcal{E} u_h  \; \mathcal{E} v_h \, ds = \\
        = \int_{\tilde{\Omega}} f v_h \, dx
        - \int_{\tilde{\Gamma}} (\nabla v_h \cdot \tilde{\mathbf{n}}) {g} \, ds
        + \int_{\tilde{\Gamma}} \sigma_\Gamma g v_h \, ds.
    \end{split}
\end{equation}
where $\tilde{\mathbf{n}}$ is the outward normal to $\tilde{\Gamma}$ and $\sigma_\Gamma$ is a penalty parameter.
We notice that our discrete solution $u_h$ is a piecewise polynomial function defined on the
background mesh $\mathcal{T}_h$, hence its Taylor expansion can be computed directly by evaluating the
function values at the points on the true boundary $\Gamma$. This allows us to avoid computation of higher-order
derivatives of $u_h$. Note that the resulting form is not symmetric due to the presence of the extrapolated
boundary condition.

\section{Multigrid Preconditioner}
\label{sec:mg_preconditioner}
The SBM formulation described above leads to a large, sparse linear system. Due to
the lack of symmetry and possible indefiniteness of the resulting matrix, we employ a Krylov subspace method,
specifically GMRES, for its solution. The convergence of GMRES is sensitive to the condition number of the system
matrix, which grows with both the number of elements in the mesh and the polynomial degree.
To accelerate the solution, we employ a multigrid preconditioner~\cite{hackbusch1985multi}. Our Cartesian
background mesh naturally facilitates the construction of a nested hierarchy of meshes
$\{\mathcal{T}_\ell\}_{\ell=0}^L$,
where $\ell$ denotes the level and $L$ is the finest level. On each mesh $\mathcal{T}_\ell$, we define a finite element
space $\mathbb{V}_\ell$ consisting of continuous piecewise polynomials of a fixed degree $p$.

This represents an important difference from the work on a Discontinuous Galerkin SBM
preconditioner~\cite{wichrowski2025geometric}, where an hp-multigrid strategy with lower polynomial degrees on coarser
levels was necessary to achieve good performance. In the present continuous Galerkin context, we find that a simpler
h-multigrid approach is enough. We note that while the background meshes are nested, the sets of active cells on
different levels (determined by the threshold $\lambda$ for their volume fraction outside the domain) are not
necessarily nested. This lack of geometric nestedness can affect the efficiency of standard inter-grid transfer
operators.

With this hierarchy of meshes and spaces established, the multigrid method is defined through three crucial components:
a smoother, which reduces high-frequency error components on each level; transfer operators, which move information
between different resolution levels; and a coarse-grid solver. For preconditioning, we use a single multigrid V-cycle.

\subsection{Smoother}
\label{sec:smoother}
We first consider Richardson iterations with preconditioners $P$. Given the current
approximation $u_h$ and the right-hand side $f_h$, a~smoothing step updates
the approximation:
\[
    u_h^{k+1} = u_h^k + P(f_h - A_h u_h^k).
\]
We decompose the space $\DGell$ into a sum of subspaces $\mathbb{W}_i$, i.e., $\DGell = \sum_{i=1}^N \mathbb{W}_i$. Let
$R_i\colon \DGell \to \mathbb{W}_i$ be the restriction operator and $A_i$ be the restriction of $A$ to the subspace
$\mathbb{W}_i$. Then, the additive subspace correction preconditioner is defined as:
\[
    B = \omega \sum_{i=1}^N R_i^T A_i^{-1} R_i
\]
where $\omega$ is a relaxation parameter. Alternatively, the successive subspace correction method is a subspace
correction method where the subspaces $\mathbb{W}_i$ are visited in a sequential manner. The procedure can be performed
by applying one Richardson iteration with the preconditioner $P_i$ for each subspace $\mathbb{W}_i$. Since in each step
only one subspace is corrected, the residual only changes locally and an efficient implementation is possible. The
preconditioner updates the solution by sequentially applying corrections for each subspace. For each subspace
$\mathbb{W}_i$, a correction is computed and applied:
\[
    u_h \leftarrow u_h + R_i^T A_i^{-1} R_i (f_h - A_h u_h).
\]
This process is repeated for all subspaces $i=1, \ldots, N$. If the subspaces are chosen as one-dimensional spaces
spanned by the basis functions, then the preconditioner is equivalent to either Jacobi (additive) or Gauss-Seidel
(successive).

Viewing the smoother as a subspace correction method provides a useful framework. The effectiveness of such a smoother
depends critically on the choice of subspaces and the local problems solved on them. For continuous finite element
methods, simple choices like one-dimensional subspaces corresponding to individual degrees of freedom (leading to
Jacobi or Gauss-Seidel smoothers) are ineffective for SBM. These methods implicitly enforce a zero correction on the
boundary of the single-vertex \emph{patch}, which conflicts with the non-standard boundary terms introduced by SBM, a
fact confirmed by our preliminary experiments.

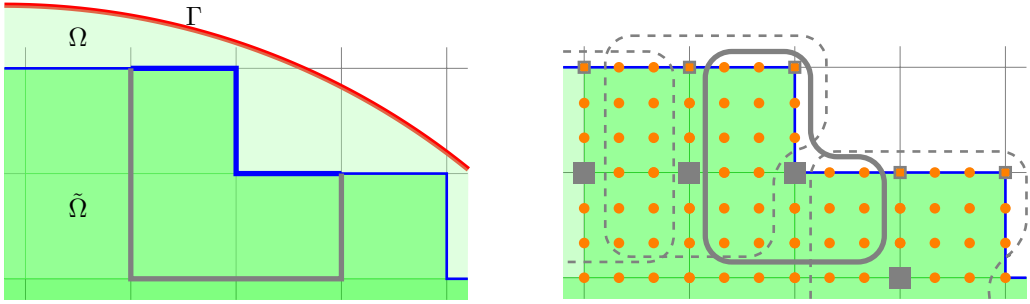
\begin{figure}[!ht]
    \centering
    \begin{subfigure}[b]{0.4\linewidth}
        \centering
\begin{tikzpicture}[scale=1.4]
    % Draw the Cartesian grid
    \draw[step=1cm, gray, very thin] (0.8,-0.2) grid (5.2,2.2);

    % Highlight cells inside the domain (interior cells)
    \foreach \x/\y in {2/1, 1/0, 1/1, 2/0, 3/0, 4/0}
        {
            \fill[green, opacity=0.5] (\x,\y) rectangle (\x+1,\y+1);
        }
    \foreach \x/\y in {1/0 }{
            \fill[green, opacity=0.5] (\x,\y) rectangle (\x-0.2,\y+2);
        }

    \fill[green, opacity=0.5] (0.8,-0.2) rectangle (5.2,0);

    \draw[line width=2pt, red, domain=51:90, samples=100, variable=\t]
    plot ({0.8 + 7.0*cos(\t)}, {-4.4 + 7.0*sin(\t)});
    \node at (2.6,2.5) {$\Gamma$};

    % Fill the area below the red line with light green to indicate the domain Omega
    \fill[green!30, opacity=0.4] plot[domain=51:90, samples=100, variable=\t] ({0.8 + 7.0*cos(\t)}, {-4.4 +
            7.0*sin(\t)}) -- (0.8,0) -- (5.2,0) -- cycle;

    % Draw the surrogate boundary (only the top boundary of the interior cells)
    \draw[line width=1pt, blue] (0.8,2) -- (3,2) -- (3,1) -- (4,1) -- (5,1)-- (5,0)-- (5.2,0.) ;

    % Draw the surrogate boundary in the single patch case
    \draw[line width=2pt, blue] (2,2) -- (3,2) -- (3,1) -- (4,1);

    % Draw the interior boundary 
    \draw[line width=2pt, gray] (2,2) -- (2,0) -- (4,0) -- (4,1);

    \node at (1.5,0.7) {$\tilde{\Omega}$};

    \node at (1.5,2.3) {${\Omega}$};

\end{tikzpicture}
    \end{subfigure}
    \qquad
    \begin{subfigure}[b]{0.4\linewidth}
        \centering
\begin{tikzpicture}[scale=1.4,
        dof/.style={circle, fill=black, inner sep=0.8pt}]
    % Draw the Cartesian grid
    \draw[step=1cm, gray, very thin] (0.8,-0.2) grid (5.2,2.2);

    \usetikzlibrary{calc}
    \newcommand{\drawCubicDoFsAt}[2]{%
        \foreach \xi in {0,0.33,0.66,1} {%
                \foreach \eta in {0,0.33,0.66,1}  {%
                        \node[dof, minimum size=4pt, inner sep=0pt, fill=orange] at ($(#1,#2) + (\xi,\eta)$)
                        {};%
                    }%
            }%
    }

    %Draw central vertex 
    %BTW: The name is not mine, GPT5 came up with it xD
    \newcommand{\drawConfidentSquareAt}[2]{%
        \fill[gray] ($(#1,#2)+(-0.1,-0.1)$) rectangle ($(#1,#2)+(0.1,0.1)$);%
    }

    \newcommand{\drawShySquareAt}[2]{% 
        \draw[gray, line width=1pt, fill=none, draw opacity=1] ($(#1,#2)+(-0.05,-0.05)$) rectangle
        ($(#1,#2)+(0.05,0.05)$);% 
    }

    \foreach \x/\y in {1/0 }{
            \fill[green, opacity=0.2] (\x,\y) rectangle (\x-0.2,\y+2);
        }

    \fill[green, opacity=0.4] (0.8,-0.2) rectangle (5.2,0);
    % Highlight cells inside the domain (interior cells)
    \foreach \x/\y in {2/1, 1/0, 1/1, 2/0, 3/0, 4/0}
        {
            \fill[green, opacity=0.4] (\x,\y) rectangle (\x+1,\y+1);
        }

    % Draw the surrogate boundary (only the top boundary of the interior cells)
    \draw[line width=1pt, blue] (0.8,2) -- (3,2) -- (3,1) -- (4,1) -- (5,1)-- (5,0)-- (5.2,0.) ;

    \foreach \x/\y in {2/1, 1/0, 1/1, 2/0, 3/0, 4/0}
        {
            \drawCubicDoFsAt{\x}{\y}
        }

    % Draw the patch boundary in the shy patches  
    \draw[line width=2pt, gray, rounded corners=10pt, line join=round] (2.15,2.15) -- (2.15,0.15) -- (3.85,0.15) --
    (3.85,1.15)--(3.15,1.15)--(3.15,2.15)--cycle;

    %  Draw shy center at (3,1)
    \drawConfidentSquareAt{3}{1}

    % Draw the patch boundary in the shy patches  
    \draw[line width=1pt, gray, dashed, rounded corners=10pt, line join=round] (3.15,-0.2) -- (3.15,1.2) --
    (5.20,1.2) --(5.20,0.3) --(4.82,0.1)--(4.85,-0.2);

    %  Draw shy center at (3,1)
    \drawConfidentSquareAt{4}{0}

    % Draw the patch boundary in the shy patches  
    \draw[line width=1pt, gray, dashed, rounded corners=10pt, line join=round] (0.8,2.15) -- (1.85,2.15) --
    (1.850,0.15) --(0.8,0.15);
    \drawConfidentSquareAt{1}{1}

    \draw[line width=1pt, gray, dashed, rounded corners=10pt, line join=round] (1.2,2.3) -- (3.3,2.3) --
    (3.3,1.3) --(2.8,1.15)--(2.8,0.2)--(1.2,0.2)--cycle;
    \drawConfidentSquareAt{2}{1}

    \drawShySquareAt{1}{2}
    \drawShySquareAt{2}{2}
    \drawShySquareAt{3}{2}
    \drawShySquareAt{4}{1}
    \drawShySquareAt{5}{1}

\end{tikzpicture}
    \end{subfigure}
    \caption{
        Left: Boundaries of a vertex patch. Interior boundaries shown in light gray, and the thick blue
        segments
        indicate the exterior faces where SBM boundary conditions are applied. Right: the same vertex patch with a gray
        frame
        highlighting which degrees of freedom form the patch subspace; three other patches are indicated with gray
        dashed
        outlines, and the central vertices are marked by filled gray squares. The vertices marked with small gray
        squares may  be too shy to form their own patches, depending on the shyness threshold. Note that as long as
        the shyness threshold is at most 3 all degrees of freedom from those patches will be include in other patches.
    }
    \label{fig:single_patch_boundaries}
\end{figure}

This suggests that the subspaces must be large enough to capture the local behavior of the SBM formulation. We
therefore construct local problems on larger subspaces, defined by patches of elements centered around a vertex, as
illustrated in Figure~\ref{fig:single_patch_boundaries}. By solving a local problem on this patch that includes the
relevant SBM boundary terms where applicable, we can compute a more effective correction. This approach aims to make
the local problem on the patch a more consistent approximation of the global problem, leading to a more robust
smoother.

\subsection{Full-Residual Shy Patches}
\label{sec:shy_patches}
The idea behind our proposed smoother is to construct local problems on patches of cells that are consistent with
the global problem, particularly in handling the non-standard terms introduced by the SBM. We start by constructing
vertex patches, where a patch consists of all active cells connected to a given vertex. The subspace for the smoother
correction is then defined by the degrees of freedom (DoFs) supported by the cells within this patch.

Next we ensure that the local problem solved on the patch is a faithful representation of the global problem. To
achieve this, we only include DoFs for which the complete residual can be evaluated using information solely from the
cells within the patch. This means we exclude any DoFs whose support extends to cells outside the current patch. By
doing so, we guarantee that the local correction is computed using the true global residual, avoiding inconsistencies
that would arise from incomplete information, especially near the surrogate boundary.
Figure~\ref{fig:single_patch_boundaries} illustrates this concept for a few patches.

We further refine this patch selection with a concept we call \emph{shyness}. Much like a shy person at a party needs a
circle of friends to feel confident enough to join the dance, our vertices are considered shy. A vertex will only form
the center of a patch if it is surrounded by a sufficient number of \emph{friends} --- in this case, active cells. We
define a shyness threshold as the minimum number of active cells required around a vertex to consider it for patch
construction. If a vertex does not meet this threshold, it is deemed too isolated, and no patch is formed around it.
This strategy prevents the creation of small local problems, which reduces the computational cost of the smoother and
avoids potential issues with ill-conditioning that can arise from patches with too few active cells.

This raises an important question about the maximum shyness threshold required to ensure that every degree of freedom
is included in at least one patch. If a DoF is not part of any patch, its value will never be updated by the smoother,
which is detrimental to convergence. In two dimensions, a vertex can be adjacent to at most four cells. As illustrated
in Figure~\ref{fig:single_patch_boundaries}, if a vertex is adjacent to only one or two active cells, the DoFs
associated with that vertex might not be included in any other patch if their respective central vertices are also
\emph{shy}. To guarantee that every DoF belongs to at least one patch, the shyness threshold $\xi$ must be at most 3 in
2D, while in 3D it must be at most 4. This ensures that even if a vertex is too shy to form its own patch, its
associated DoFs are guaranteed to be included in the patch of a sufficiently \emph{sociable} neighboring vertex.

% While it is theoretically possible
% for some DoFs to be left out with a lower threshold, which would hinder convergence, setting the minimum shyness to 3
% provides a robust covering of all DoFs.

In~\cite{wichrowski2025geometric}, satisfactory results were obtained with 3 smoothing steps. A multi-stage smoothing
strategy was proposed in~\cite{cui2025multigrid} for CutFEM, where an initial global smoothing pass was followed by
additional passes restricted to patches containing cut cells. We adopt a similar strategy here. The first smoothing
step is applied to all patches in the domain. Subsequent steps are then selectively applied only to those patches that
contain cells adjacent to the surrogate boundary. This focuses the computational effort of the smoother on the region
where the SBM introduces non-standard terms and where the error is often most difficult to resolve. These
boundary-adjacent patches are identified once during the construction of the smoother.

\subsection{Transfer Operators}
\label{sec:transfer_operators}
The transfer of information between different levels of the multigrid hierarchy is handled by prolongation and
restriction operators. The prolongation operator, $I_{\ell-1}^\ell$, maps a function from the coarse space
$\mathbb{V}_{\ell-1}$ to the fine space $\mathbb{V}_\ell$, while the restriction operator, $I_\ell^{\ell-1}$, transfers
a function from the fine space to the coarse space. For standard multigrid methods on nested meshes with continuous
finite elements, prolongation is typically the natural embedding of the coarse function space into the fine one, and
restriction is its transpose.

A significant challenge in applying multigrid to unfitted methods like SBM is the lack of geometric nestedness of the
computational domains. The set of active cells $\tilde{\Omega}_\ell$ on a given level $\ell$ is determined
independently based on the intersection with the true domain $\Omega$. Consequently, the active domain on the fine
level, $\tilde{\Omega}_L$, is not necessarily a subset of the active domain on a coarser level, $\tilde{\Omega}_{\ell}$
(when viewed on the fine grid). This can lead to inconsistencies where, for example, an active fine-grid cell
corresponds to a non-active coarse-grid cell.

This lack of nestedness can impair the effectiveness of standard transfer operators~\cite{wichrowski2025geometric}.
Information from an active region on the fine grid might be restricted to a non-active region on the coarse grid, where
it is essentially discarded, breaking the flow of information required for an efficient multigrid cycle.

\section{Implementation details}
\label{sec:implementation}
The numerical implementation of our multigrid solver is based on the open-source finite element library
\texttt{deal.II}~\cite{dealii2024}. It provides a comprehensive framework for the implementation of~finite
element methods, including mesh handling, finite element spaces, assembly of linear systems, and interfaces to various
linear algebra solvers and preconditioners. In this paper we build upon the implementation of the multigrid solver for
DG-SBM presented in~\cite{wichrowski2025geometric}.

\subsection{Background mesh and geometry handling}
When implementing SBM on a non-body-fitted mesh, we need to handle cells intersected by the true boundary $\Gamma$. We
use a level set function $\phi(\mathbf{x})$ to implicitly define the domain $\Omega = \{\mathbf{x} \mid
    \phi(\mathbf{x}) < 0\}$, with $\Gamma = \{\mathbf{x} \mid \phi(\mathbf{x}) = 0\}$. Cells of the background mesh are
classified based on their intersection with the zero level set: cells entirely inside $\Omega$ (interior), cells
entirely outside $\Omega$ (exterior), and cells intersected by $\Gamma$. Then, for the intersected cells the fraction
of the cell volume inside $\Omega$ is computed, and the cell is classified as active if this fraction is greater than a
threshold $\lambda$. This is handled using non-matching quadrature rules implemented in \texttt{deal.II}, which are
based on the techniques described in~\cite{saye2015high}.

In our approach, degrees of freedom are formally assigned to all cells of the background mesh, including those that are
classified as non-active. In the matrix assembly we ignore the contributions from these non-active cells, effectively
removing them from the system. This results in a singular global matrix; however, since the corresponding degrees of
freedom do not influence the solution in the active domain, this does not pose a problem for the iterative solver. The
only affected part of the multigrid is the coarse-grid solver, which we handle by using a direct solver. To make the
coarse system invertible, we set one on the diagonal entries of zero rows corresponding to non-active DoFs.

\subsection{Processing surrogate boundary and matrix assembly}
\label{sec:sbm_assembly}
The SBM requires computing the closest point projection from points on the surrogate boundary $\tilde{\Gamma}$ to the
true boundary $\Gamma$. This projection is found by solving a local nonlinear optimization problem for each quadrature
point on $\tilde{\Gamma}$, which minimizes the distance to the true boundary defined by the zero level set of a
function $\phi(\mathbf{x})$. The full details of the formulation, which involves a Lagrange multiplier approach solved
with a Newton-Raphson method, are described in~\cite{wichrowski2025geometric}. To ensure sufficient
smoothness for the derivative calculations required by the solver, the level set function $\phi(\mathbf{x})$ is
represented using finite elements of degree 2 for $p=1$ and degree $p$ for $p>1$. The search for the closest point is
performed only in the interior of the cells adjacent to the surrogate
boundary. While this does not guarantee that the closest point is found inside the cell, it is expected that
even if the closest point is outside the cell, the extrapolation will still yield a good approximation of the boundary
condition.

The extrapolation of the function values from the true boundary to the surrogate boundary required for the matrix
assembly process is accomplished by evaluating the values of the basis functions at the points on the surrogate
boundary. The assembly of the cell contributions and interior faces to the system matrix and right-hand side vector is
performed using the standard finite element assembly process provided by \texttt{deal.II}.

\subsection{Multigrid structures}

For the multigrid hierarchy, we leverage \texttt{deal.II}'s built-in capabilities for handling nested meshes and
defining finite element spaces on each level. The standard projection operators provided by \texttt{deal.II} are used
for the prolongation ($I_{\ell-1}^\ell$) and restriction ($I_\ell^{\ell-1}$) operators, transferring data between
coarser and finer grid levels. These operators act on the entire background mesh, transferring the solution for all
degrees of freedom, irrespective of whether they correspond to active or non-active cells. This involves transferring
values in regions outside the computational domain. Furthermore, since the smoother and residual evaluations are
restricted to the active cells, the values in the inactive regions do not propagate into the solution within the domain
of interest, nor do they influence the convergence of the method.

The assembly of the system matrix $A_\ell$ on each level $\ell$ of the multigrid hierarchy follows a procedure
analogous to that on the finest level. The level set function defining the domain geometry is interpolated onto the
mesh $\mathcal{T}_\ell$. Based on this interpolated level set and the chosen threshold $\lambda$, active cells for
level $\ell$ are identified. The SBM bilinear form is then used to assemble the local contributions to $A_\ell$ only
for these active cells. The cells deemed non-active on level $\ell$ are ignored during the assembly, effectively
decoupling their degrees of freedom from the system on that level.

\section{Numerical results}
\label{sec:numerical_results}
In this section, we present numerical results to evaluate the performance of the proposed SBM multigrid
preconditioner for the Poisson equation. We investigate its effectiveness in terms of convergence rates, iteration
counts under mesh refinement ($h$-refinement) and polynomial degree increase ($p$-refinement). To better illustrate the
multigrid preconditioner's performance characteristics and enable meaningful comparison of iteration counts, we solve
the system with a tolerance of $10^{-8}$ for the relative residual reduction in all experiments.
Solver failure is defined as exceeding 100 GMRES iterations without reaching this tolerance.
The initial background mesh is a square (or cube in 3D) domain covering the range $[-1.01, 1.01]$ in each dimension,
subdivided into $4$ cells in each coordinate direction at refinement level $L=0$. We use a penalty parameter $\sigma = 5$ in the SBM
formulation. This parameter is selected to be sufficiently large to ensure stability of the Nitsche coupling on the
surrogate boundary.

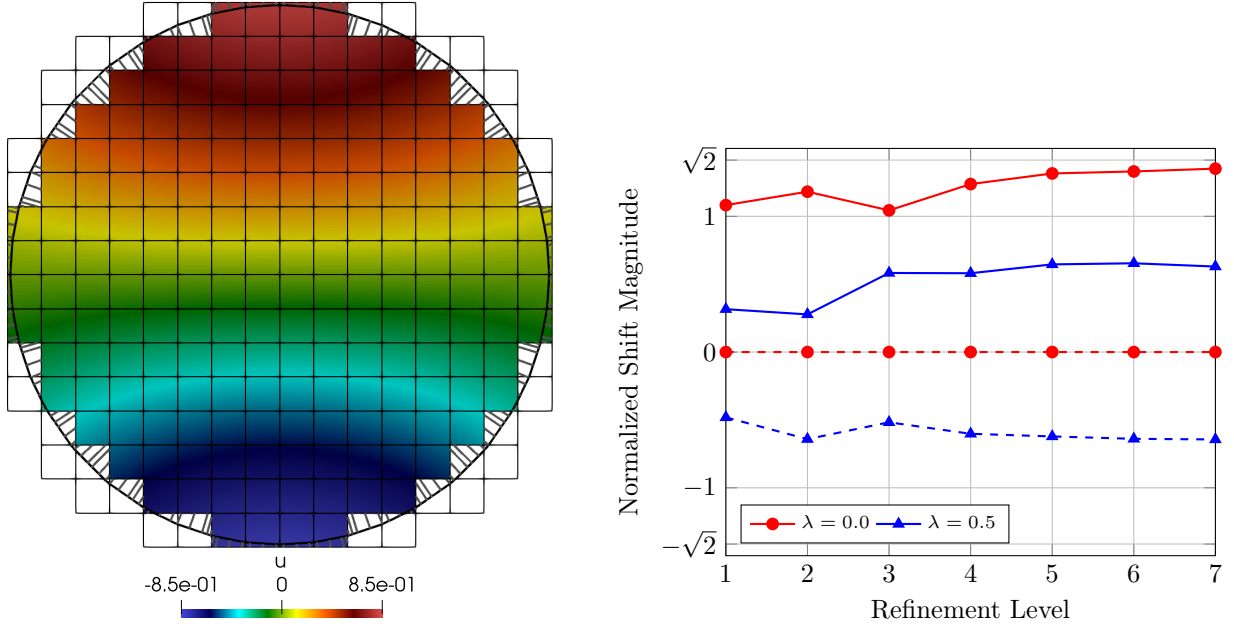
\begin{figure}[h!]
    \centering
    % Placeholder for shift distribution plot
    \begin{subfigure}{0.49\textwidth}
        \includegraphics[width=\textwidth]{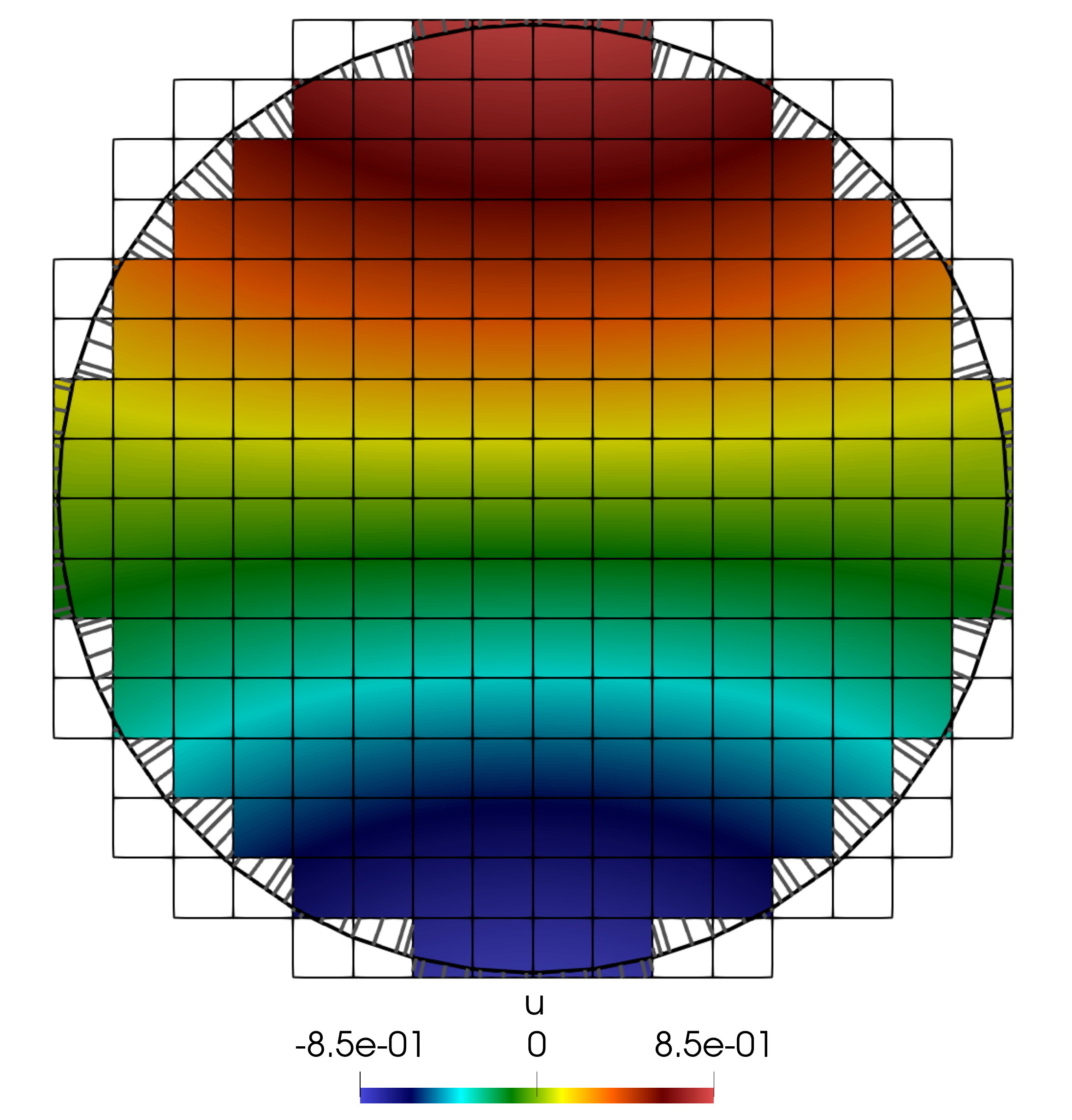}
    \end{subfigure}
    \hfill
    \begin{subfigure}{0.49\textwidth}
        \centering
        %2D
        \begin{tikzpicture}
            \begin{axis}[
                    xlabel={Refinement Level},
                    ylabel={Normalized Shift Magnitude},
                    xmin=1, xmax=7,
                    xtick={1,2,3,4,5,6,7},
                    ymin=-1.5,
                    ymax=1.5,
                    ytick={-1.41421356, -1, 0,1, 1.41421356},
                    yticklabels={$-\sqrt{2}$,$-1$ ,$0$,  $1$, $\sqrt{2}$},
                    width=\textwidth, % Use \textwidth for subfigure's width
                    legend pos=south west, % Adjusted legend position
                    legend style={font=\scriptsize, cells={anchor=west}, legend columns=2},
                    height=7cm, % Adjust height as needed
                    grid=major,
                ]
                % Data from 2D_hp_stand_ball_p3.txt
                % Lambda = 0.0
                \addplot+[mark=*,color=red, mark options={solid}, thick,forget plot,  dashed ] coordinates {
                        (1, 0.0000) (2, 0.0000) (3, 0.0000) (4, 0.0000) (5, 0.0000) (6, 0.0000) (7, 0.0000)
                    };
                \addlegendentry{$\lambda=0.0$}
                \addplot+[mark=*,color=red, mark options={solid}, thick ] coordinates {
                        (1, 1.0825) (2, 1.1813) (3, 1.0437) (4, 1.2373) (5, 1.3157) (6, 1.3304) (7, 1.3510)
                    };

                % Lambda = 0.25
                % \addplot+[mark=square*,color=darkgreen, mark options={solid}, thick, forget plot,  dashed ] coordinates
                %     {
                %         (1, -0.2334) (2, -0.3089) (3, -0.3904) (4, -0.3152) (5, -0.4135) (6, -0.4265) (7, -0.4302)
                %     };
                % \addlegendentry{$\lambda=0.25$}
                % \addplot+[mark=square*,color=darkgreen, mark options={solid}, thick] coordinates {
                %         (1, 0.7755) (2, 0.8005) (3, 0.7998) (4, 0.8367) (5, 0.8399) (6, 0.8450) (7, 0.8609)
                %     };

                % Lambda = 0.5
                \addplot+[mark=triangle*,color=blue, mark options={solid}, thick, forget plot,	dashed ] coordinates {
                        (1, -0.4811) (2, -0.6414) (3, -0.5174) (4, -0.6021) (5, -0.6216) (6, -0.6392) (7, -0.6436)
                    };
                \addlegendentry{$\lambda=0.5$}
                \addplot+[mark=triangle*,color=blue, mark options={solid}, thick] coordinates {
                        (1, 0.3159) (2, 0.2774) (3, 0.5826) (4, 0.5807) (5, 0.6458) (6, 0.6540) (7, 0.6299)
                    };

            \end{axis}
        \end{tikzpicture}
    \end{subfigure}
    \caption{Left: Exemplary numerical solution for the Poisson equation on a unit ball, visualized on the
        surrogate domain ($\lambda=0.25$), obtained with the shifted boundary method. The true boundary $\Gamma$ is
        shown in black lines; connections between quadrature points ($p=2$) on the surrogate boundary and their
        projections on the true boundary are illustrated with grey lines. Note that with $\lambda=0.25$,
        some parts of the true boundary lie within the surrogate domain.
        Right: Illustrative distribution of minimum and maximum shift magnitudes (normalized by cell size $h$) on
        the surrogate boundary for the unit ball problem with $p=3$. The values $\pm\sqrt{2}$ represent the
        theoretical maximum possible normalized shift magnitude in 2D for a square cell. For $\lambda=0.0$ the maximum
        shift magnitude was 1.3510, while for $\lambda=0.5$ it was 0.6436 and the minimum -0.6299.
    }
    \label{fig:shift_distribution}
\end{figure}

The majority of our tests are conducted on a unit ball $\Omega = \{\mathbf{x} \in \mathbb{R}^d : \|\mathbf{x}\| < 1\}$,
which corresponds to a unit disk in 2D and a unit ball in 3D (the 2D case is depicted on the left panel of
Fig.~\ref{fig:shift_distribution}). We solve the Poisson equation \eqref{eq:poisson} with a constant right-hand side
$f=1$ and homogeneous Dirichlet boundary conditions $u=0$ on $\Gamma = \partial\Omega$. While geometrically simple, the
unit ball serves as an insightful benchmark. Any sufficiently smooth ($C^1$) complex boundary, when viewed at a fine
enough mesh resolution, locally resembles a flat plane. The unit ball, due to its uniform curvature, presents a
comprehensive range of intersection angles between the true boundary and the background mesh cells. Furthermore, the
boundary intersects cells at various locations relative to cell centers and faces, leading to a diverse distribution of
shift vector magnitudes and directions. This includes scenarios where the shift vector points from the surrogate
boundary towards the interior of the true domain (which we denote as negative shifts if they oppose the outward normal
of the surrogate boundary). The right panel of Figure~\ref{fig:shift_distribution} depicts the minimum and maximum
shift magnitudes observed across the surrogate boundary for a typical discretization. The shift magnitudes are
normalized by the cell size $h$; for a unit ball, the theoretically largest possible shift magnitude in 2D is
$\sqrt{2}h$, occurring if the surrogate boundary point is at a cell corner and the true boundary passes through the
diagonally opposite corner.

The results were computed on a machine with AMD Epyc 7282 with 256GB of RAM. Due to prohibitive memory requirements of
high order discretizations in 3D, some of the larger problems could not be run to completion.

The primary metric for evaluating the multigrid preconditioner is the number of GMRES iterations required to reduce the
initial residual by a factor of $10^{-8}$. If the threshold is not met within 100 iterations, the solver is considered
to have failed. All experiments are conducted using the implementation described in the previous section. As our solver
depends on the value of the threshold parameter $\lambda$, we will first explore the solver performance on 2D problems
as they are less computationally intensive and then show some results for 3D problems with tuned $\lambda$.

\subsection{Iteration counts}

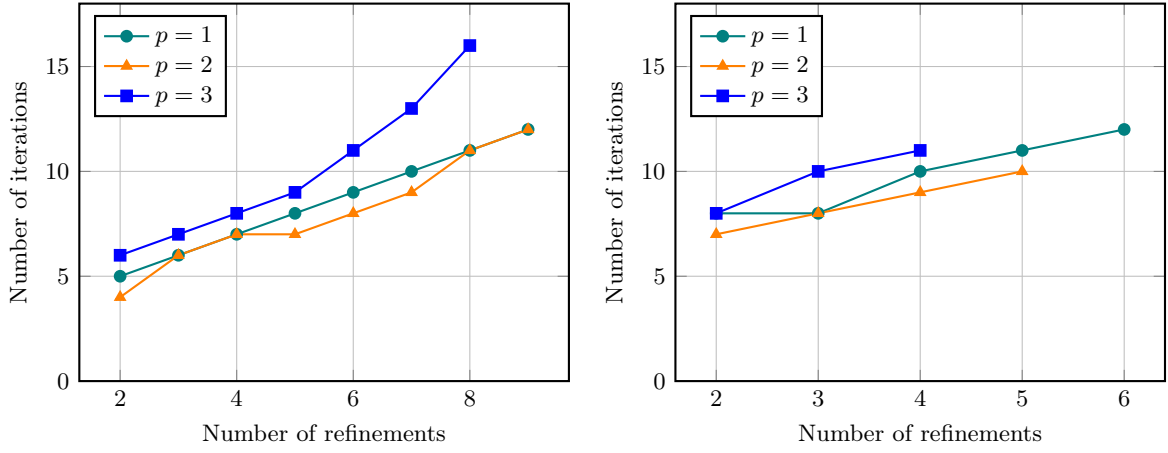
\begin{figure}[htbp]
    \centering
    \begin{tikzpicture}
    \begin{scope}
        \begin{axis}[
                width=0.45\columnwidth,
                xlabel={Number of refinements},
                ylabel={Number of iterations},
                ymin=0,
                ymax =18,
                legend pos=north west,
                paperplot]

            \immediate\write18{python3 scripts/compute_iterations.py results/full_2D.csv 3 3. 0.0 > results/tmp/iterations_2D.csv}

            \newcommand{\addIterPlot}[1]{%
                \addplot+[thick] table[
                        col sep=comma,
                        x=refinements,
                        y=iterations,
                        restrict expr to domain={\thisrow{p}}{#1:#1},
                        unbounded coords=discard
                    ] {results/tmp/iterations_2D.csv}; \addlegendentry{$p=#1$}%
            }
            \addIterPlot{1}
            \addIterPlot{2}
            \addIterPlot{3}

        \end{axis}
    \end{scope}

    \begin{scope}[xshift=0.48\columnwidth]
        \begin{axis}[
                width=0.45\columnwidth,
                xlabel={Number of refinements},
                ylabel={Number of iterations},
                ymin=0,
                ymax =18,
                legend pos=north west,
                paperplot]

            \immediate\write18{python3 scripts/compute_iterations.py results/full_3D.csv 3 3. 0.0 > results/tmp/iterations_3D.csv}

            \newcommand{\addIterPlotThree}[1]{%
                \addplot+[thick] table[
                        col sep=comma,
                        x=refinements,
                        y=iterations,
                        restrict expr to domain={\thisrow{p}}{#1:#1},
                        unbounded coords=discard
                    ] {results/tmp/iterations_3D.csv}; \addlegendentry{$p=#1$}%
            }
            \addIterPlotThree{1}
            \addIterPlotThree{2}
            \addIterPlotThree{3}

        \end{axis}
    \end{scope}
\end{tikzpicture}
    \caption{Iteration counts for multigrid preconditioned GMRES solver using the Full-Residual Shy Patch smoother with
        varying shyness thresholds equal to 3 on unit ball problem for polynomial degrees $p=1,2,3$. Number of smoothing steps is fixed to 3. Left panel presents results in 2D with shyness $3$, right panel in 3D with shyness $4$. In 3D some data points are missing due to excessive memory requirements.}
    \label{fig:iteration_counts_shy3_smoothing3}
\end{figure}

Figure~\ref{fig:iteration_counts_shy3_smoothing3} illustrates the iteration counts with respect to the number of
refinement levels for both 2D and 3D problems. The results demonstrate the $h$-robustness of the proposed
preconditioner, as the iteration counts remain largely independent of the mesh size for a fixed polynomial degree. We
observe a mild increase in iterations as the polynomial degree $p$ increases, which is expected when using
$h$-multigrid for high-order discretizations without $p$-coarsening. However, patch-smoother for matching meshes has
been shown to be robust even in $p$~\cite{wichrowski2025smoothers,pavarino1993additive}, suggesting that further
improvements may be possible with additional refinements to the smoother or the multigrid strategy.

We further investigate the influence of the number of smoothing steps $s$ on the solver performance.
Table~\ref{tab:iter_smoothing} presents the iteration counts for varying $s$ with a fixed shyness threshold $\xi=3$. In
the 2D case (Table~\ref{tab:iter_2d_smooth}), for $p=1$, the solver is robust even with a single smoothing step.
However, for $p=2$ and $p=3$, increasing the number of smoothing steps is crucial. For $p=3$, the solver fails with
$s=1$ but converges rapidly with $s=2$ or $s=3$. A similar trend is observed in the 3D case
(Table~\ref{tab:iter_3d_smooth}), where increasing $s$ consistently reduces the iteration count, particularly for
higher polynomial degrees.

\begin{table}[htbp]
    \centering
    \caption{Iteration counts for varying smoothing steps $s$ with fixed shyness threshold $\xi=3$. Left: 2D problem ($L=8$). Right: 3D problem ($L=4$).}
    \label{tab:iter_smoothing}

    \begin{subtable}{0.45\textwidth}
        \centering
        \caption{2D problem ($L=8$)}
        \label{tab:iter_2d_smooth}
        \begin{tabular}{c|ccc}
            \toprule
            $p \setminus s$ & 1   & 2  & 3  \\
            \midrule
            1               & 11  & 11 & 11 \\
            2               & 22  & 13 & 11 \\
            3               & --- & 29 & 16 \\
            \bottomrule
        \end{tabular}
    \end{subtable}
    \begin{subtable}{0.45\textwidth}
        \centering
        \caption{3D problem ($L=4$)}
        \label{tab:iter_3d_smooth}
        \begin{tabular}{c|ccc}
            \toprule
            $p \setminus s$ & 1  & 2  & 3  \\
            \midrule
            1               & 13 & 11 & 10 \\
            2               & 15 & 10 & 9  \\
            3               & 28 & 15 & 11 \\
            \bottomrule
        \end{tabular}
    \end{subtable}
\end{table}

Finally, we analyze the effect of the shyness threshold $\xi$ on the convergence, with the number of smoothing steps
fixed at $s=3$. The results are summarized in Table~\ref{tab:iter_shyness}. For the 2D problem
(Table~\ref{tab:iter_2d_shy}), we present results for the baseline threshold $\xi=3$. In the 3D case
(Table~\ref{tab:iter_3d_shy}), we compare thresholds $\xi=3, 4, 5$. The results indicate that while the method is
robust for $\xi=3$ and $\xi=4$, increasing the threshold to $\xi=5$ results in solver failure. This suggests that a
threshold of $\xi=5$ is too restrictive, leaving some degrees of freedom outside of any patch, thereby degrading the
quality of the smoother.

\begin{table}[htbp]
    \centering
    \caption{Iteration counts for varying shyness threshold $\xi$ with fixed smoothing steps $s=3$. Left: 2D problem ($L=8$). Right: 3D problem ($L=4$).}
    \label{tab:iter_shyness}

    \begin{subtable}{0.45\textwidth}
        \centering
        \caption{2D problem ($L=8$)}
        \label{tab:iter_2d_shy}
        \begin{tabular}{c|cccc}
            \toprule
            $p \setminus \xi$ & 1  & 2  & 3  & 4   \\
            \midrule
            1                 & 11 & 11 & 11 & --- \\
            2                 & 10 & 10 & 10 & --- \\
            3                 & 15 & 15 & 15 & --- \\
            \bottomrule
        \end{tabular}
    \end{subtable}
    \begin{subtable}{0.45\textwidth}
        \centering
        \caption{3D problem ($L=4$)}
        \label{tab:iter_3d_shy}
        \begin{tabular}{c|ccc}
            \toprule
            $p \setminus \xi$ & 3  & 4  & 5   \\
            \midrule
            1                 & 10 & 10 & --- \\
            2                 & 9  & 9  & --- \\
            3                 & 11 & 11 & --- \\
            \bottomrule
        \end{tabular}
    \end{subtable}
\end{table}

We also examine the influence of the cell threshold parameter $\lambda$ on the solver performance.
Table~\ref{tab:iter_lambda} presents the iteration counts for varying $\lambda$ with fixed shyness threshold $\xi=3$
and smoothing steps $s=3$. In the 2D case (Table~\ref{tab:iter_2d_lambda}), for $p=1$ and $p=2$, the solver is
relatively robust to $\lambda$, although $\lambda=0.0$ yields slightly lower iteration counts. However, for $p=3$,
large values of $\lambda$ ($\lambda \ge 0.5$) lead to increased iterations or failure, while $\lambda \le 0.25$
restores convergence. In the 3D case (Table~\ref{tab:iter_3d_lambda}), the sensitivity is more pronounced. For $p \ge
    2$, $\lambda=0.5$ results in solver failure, while $\lambda=0.0$ provides robust convergence. This might be due to
imaginary eigenvalues~\cite{wichrowski2025geometric} appearing when shifts are towards the interior of the domain.

\begin{table}[htbp]
    \centering
    \caption{Iteration counts for varying cell threshold $\lambda$ with fixed shyness threshold $\xi=3$ and smoothing steps $s=3$. Left: 2D problem ($L=8$). Right: 3D problem ($L=4$).}
    \label{tab:iter_lambda}

    \begin{subtable}{0.45\textwidth}
        \centering
        \caption{2D problem ($L=8$)}
        \label{tab:iter_2d_lambda}
        \begin{tabular}{c|cccc}
            \toprule
            $p \setminus \lambda$ & 0.0 & 0.25 & 0.5 & 0.75 \\
            \midrule
            1                     & 11  & 13   & 13  & 13   \\
            2                     & 11  & 11   & 11  & 12   \\
            3                     & 16  & 14   & 27  & ---  \\
            \bottomrule
        \end{tabular}
    \end{subtable}
    \begin{subtable}{0.45\textwidth}
        \centering
        \caption{3D problem ($L=4$)}
        \label{tab:iter_3d_lambda}
        \begin{tabular}{c|cc}
            \toprule
            $p \setminus \lambda$ & 0.0 & 0.5 \\
            \midrule
            1                     & 10  & 13  \\
            2                     & 9   & --- \\
            3                     & 11  & --- \\
            \bottomrule
        \end{tabular}
    \end{subtable}
\end{table}

\subsection{p-Multigrid}

We observed a mild increase in iteration counts with increasing polynomial degree for the h-multigrid method. This is
likely due to the properties of the transfer operators between the grid levels. To obtain a method that is robust with
respect to the polynomial degree and does not rely on a hierarchy of nested meshes, we implemented a p-multigrid
strategy. In this approach, the grid is fixed, and the multigrid hierarchy is constructed by coarsening the polynomial
degree of the finite element space. The coarsest level corresponds to $p=1$, which is solved using a direct solver.
Currently, this strategy is limited to relatively coarse grids, as the direct solver for the $p=1$ system on the fine
grid becomes prohibitively expensive in 3D. However, this limitation can be overcome by using a geometric or even
algebraic multigrid solver for the coarse $p=1$ problem, which would extend the applicability of the p-multigrid method
to 3D.

Table~\ref{tab:pmg_iterations} presents the iteration counts for the p-multigrid solver for polynomial degrees $p=2, 3,
    4$ and varying mesh refinement levels. The results demonstrate excellent robustness with respect to both the mesh size
$h$ and the polynomial degree $p$. The iteration counts are very low and remain nearly constant as the mesh is refined,
confirming the effectiveness of the p-multigrid approach combined with the Shy Patch smoother. In 2D for $p=5$, the
solver starts to struggle. In 3D, the method remains robust for $p=2$ and $p=3$, but fails to converge for $p=4$ on the
tested meshes, suggesting that additional improvements are needed to handle higher polynomial degrees in three
dimensions.

\begin{table}[htbp]
    \centering
    \caption{Iteration counts for p-multigrid preconditioned GMRES solver with varying polynomial degree $p$ and mesh refinements. Shyness threshold $\xi=3$, smoothing steps $s=3$, cell threshold $\lambda=0.0$. Left: 2D results. Right: 3D results (from pmg\_3D.csv, gmg\_its). Dash --- indicate no converge within 100 iterations, empty entries indicate that the data is not available.}
    \label{tab:pmg_iterations}
    \begin{subtable}{0.55\textwidth}
        \centering
        \caption{2D ($L$ = refinement level)}
        \label{tab:pmg_iterations_2d}
        \begin{tabular}{c|ccccccc}
            \toprule
            $p \setminus L$ & 2  & 3  & 4  & 5  & 6  & 7  & 8 \\
            \midrule
            2               & 4  & 4  & 4  & 4  & 4  & 3  & 3 \\
            3               & 4  & 4  & 4  & 4  & 4  & 3  & 3 \\
            4               & 4  & 4  & 5  & 5  & 5  & 5  & 5 \\
            5               & 32 & 32 & 33 & 32 & 34 & 34 &   \\
            \bottomrule
        \end{tabular}
    \end{subtable}
    \hfill
    \begin{subtable}{0.40\textwidth}
        \centering
        \caption{3D ($L$ = refinement level)}
        \label{tab:pmg_iterations_3d}
        \begin{tabular}{c|ccc}
            \toprule
            $p \setminus L$ & 2   & 3   & 4   \\
            \midrule
            2               & 5   & 6   & 5   \\
            3               & 5   & 5   & 6   \\
            4               & --- & --- & --- \\
            \bottomrule
        \end{tabular}
    \end{subtable}
\end{table}

\subsection{Comparison with other methods}

\begin{table}[h!]
    \centering
    \caption{GMRES iteration counts (tolerance $10^{-8}$) for the Poisson problem. Comparison of the continuous SBM ($\lambda=0$) preconditioned with the Full-Residual Shy Patch smoother (shyness $\xi=3$, $s=3$, 2D) against DG-SBM with a cellwise SSOR smoother ($s=3$) and algebraic multigrid (AMG). Solver failures (no convergence within 100 iterations) are indicated by ``---''. Empty entries indicate that the data is not available.}
    %DG-SBM results for lambda 0.25: choice of lambda because it gives better results than 0.0
    \label{tab:comparison_cutfem}
    \begin{tabular}{lll|ccccccc}
        \hline
        Formulation                             & Preconditioner                                               & Degree & \multicolumn{7}{c}{Mesh refinement}                              \\
                                                &                                                              &        & 2                                   & 3 & 4  & 5  & 6  & 7  & 8  \\
        \hline
        %data source: results/full_2D.csv
        \multirow{7}{*}{SBM, $\lambda=0$}       & \multirow{3}{*}{h-MG,Shy-Patch smoother}                     & $1$    & 5                                   & 5 & 7  & 8  & 9  & 10 & 11 \\
                                                &                                                              & $2$    & 4                                   & 6 & 7  & 7  & 8  & 9  & 10 \\
                                                &                                                              & $3$    & 5                                   & 7 & 8  & 9  & 11 & 13 & 15 \\
        %data source: results/pmg_2d.csv
        \cline{2-10}
                                                & \multirow{3}{*}{p-MG, Shy-Patch smoother}                    & $2$    & 4                                   & 4 & 4  & 4  & 4  & 3  & 3  \\
                                                &                                                              & $3$    & 4                                   & 4 & 4  & 4  & 4  & 3  & 3  \\
                                                &                                                              & $4$    & 4                                   & 4 & 5  & 5  & 5  & 5  & 5  \\
        \cline{2-10}
                                                & \multirow{1}{*}{AMG}                                         & $1$    & 2                                   & 2 & 12 & 13 & 14 & 15 & 17 \\
        \cline{1-10}
        \multirow{2}{*}{DG-SBM, $\lambda=0.25$} & {cellwise-SSOR,}                                             & $1$    & 5                                   & 6 & 7  & 8  & 9  & 9  & 11 \\
                                                & $hp$-multigrid                                               & $2$    & 7                                   & 9 & 10 & 12 & 14 & 14 & 16 \\
        \cline{1-10}
        \multirow{3}{*}{cutFEM}                 & \multirow{3}{*}{h-MG, Patch smoother\cite{cui2025multigrid}}

                                                & $1$                                                          &        &                                     &   & 6  & 6  & 6  & 5       \\
                                                &                                                              & $2$    &                                     &   &    & 9  & 8  & 7  & 7  \\
                                                &                                                              & $3$    &                                     &   &    & 17 & 14 & 13 & 13 \\
        \hline
    \end{tabular}
\end{table}

To place the performance of the proposed SBM preconditioners in context, we compare them against other established
methods for the 2D Poisson problem. Table~\ref{tab:comparison_cutfem} summarizes the iteration counts for the
continuous SBM formulation alongside results for Discontinuous Galerkin SBM (DG-SBM), Algebraic Multigrid (AMG), and
CutFEM. We focus on the 2D case to enable a comparison across a broad range of refinement levels.

We first consider the DG-SBM formulation solved using an $hp$-multigrid method with a cellwise SSOR smoother. The
proposed SBM with the Shy-Patch smoother consistently requires fewer iterations than DG-SBM. This advantage becomes
increasingly apparent as the mesh is refined, underscoring the superior smoothing capabilities of the patch-based
approach compared to cellwise smoothers within the SBM framework.

Next, we benchmark against a standard AMG preconditioner for linear finite elements ($p=1$). While AMG is highly
efficient on coarser meshes, its iteration counts tend to grow steadily with mesh refinement. In contrast, the
geometric multigrid with the Shy-Patch smoother exhibits a much more moderate increase in iterations, suggesting better
scalability properties on fine meshes.

Finally, we compare our approach with cutFEM utilizing a patch smoother~\cite{cui2025multigrid}. Our method yields
iteration counts comparable to cutFEM. While cutFEM shows remarkable stability for $p=1$, our h-multigrid approach
experiences a slight increase in iterations. However, the p-multigrid variant (p-MG) demonstrates exceptional
robustness, outperforming all other tested methods with very low and nearly constant iteration counts across all
refinement levels and polynomial degrees.

\subsection{Timing}
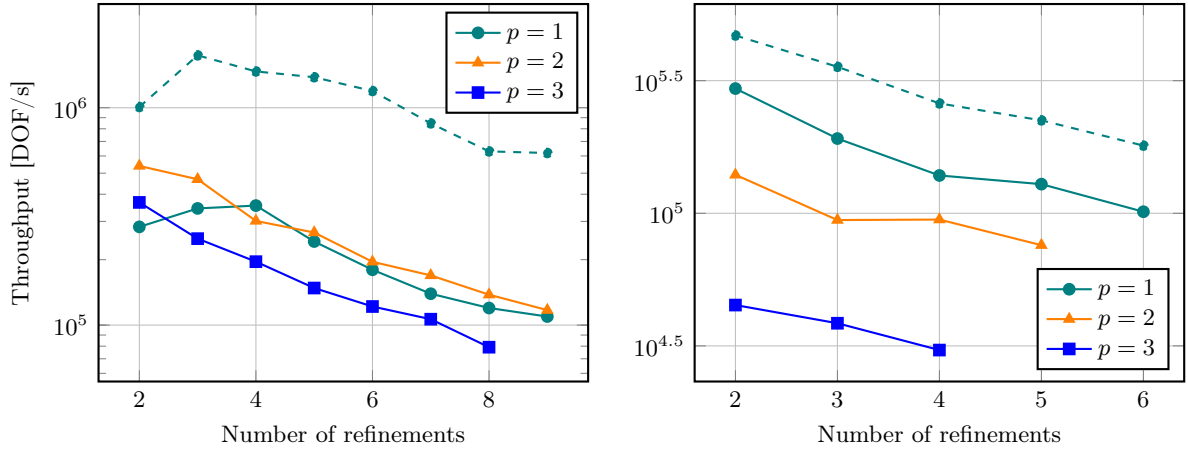
\begin{figure}[htbp]
    \centering
    \begin{tikzpicture}
    \begin{scope}
        \begin{axis}[
                width=0.45\columnwidth,
                xlabel={Number of refinements}, ylabel={Throughput [DOF/s]}, ymode=log, paperplot,
                ymax =3e6 ]

            \immediate\write18{python3 scripts/compute_throughput.py results/full_2D.csv 3 3. 0.0 > results/tmp/throughput_2D.csv}

            \newcommand{\addGMGplot}[1]{%
                \addplot+[thick] table[
                        col sep=comma,
                        x=refinements,
                        y=gmg_throughput,
                        restrict expr to domain={\thisrow{p}}{#1:#1},
                        unbounded coords=discard
                    ] {results/tmp/throughput_2D.csv}; \addlegendentry{$p=#1$}%
            }
            \addGMGplot{1}
            \addGMGplot{2}
            \addGMGplot{3}

            % AMG plots are dashed via the `dashed` plot option below.
            % Reset cycle index so AMG uses the same colors (but dashed)
            % as the GMG curves which were plotted above.
            \pgfplotsset{cycle list shift=-3}

            \newcommand{\addAMGplot}[1]{%
                \addplot+[thick, dashed, mark=*, mark options={scale=0.8}, mark repeat=1] table[
                        col sep=comma,
                        x=refinements,
                        y=amg_throughput,
                        restrict expr to domain={\thisrow{p}}{#1:#1},
                        unbounded coords=discard
                    ] {results/tmp/throughput_2D.csv};
            }
            \addAMGplot{1}

            % Restore cycle-list shift
            \pgfplotsset{cycle list shift=0}

        \end{axis}
    \end{scope}

    \begin{scope}[xshift=0.48\columnwidth]
        \begin{axis}[
                width=0.45\columnwidth,
                xlabel={Number of refinements},
                ymode=log,
                legend pos=south east,
                paperplot]
            \immediate\write18{python3 scripts/compute_throughput.py results/full_3D.csv 3 3. 0.0 > results/tmp/throughput_3D.csv}

            \newcommand{\addGMGplotThree}[1]{%
                \addplot+[thick] table[
                        col sep=comma,
                        x=refinements,
                        y=gmg_throughput,
                        restrict expr to domain={\thisrow{p}}{#1:#1},
                        unbounded coords=discard ] {results/tmp/throughput_3D.csv}; \addlegendentry{$p=#1$}%
            }

            % Add GMG curves first (match 2D pattern)
            \addGMGplotThree{1}
            \addGMGplotThree{2}
            \addGMGplotThree{3}

            % AMG plots dashed but reuse GMG colors via cycle list shift
            \pgfplotsset{cycle list shift=-3}

            \newcommand{\addAMGplotThree}[1]{%
                \addplot+[thick, dashed, mark=*, mark options={scale=0.8}, mark repeat=1] table[
                        col sep=comma,
                        x=refinements,
                        y=amg_throughput,
                        restrict expr to domain={\thisrow{p}}{#1:#1},
                        unbounded coords=discard
                    ] {results/tmp/throughput_3D.csv};
            }
            \addAMGplotThree{1}

            % Restore cycle-list shift
            \pgfplotsset{cycle list shift=0}

        \end{axis}
    \end{scope}
\end{tikzpicture}
    \caption{Comparison of throughput (DoF/s) for the proposed multigrid preconditioner with Full-Residual Shy Patch smoother with  $\xi=3$, $s=3$ (solid lines), against
        an AMG preconditioner (dashed line) from \texttt{deal.II}'s interface to \texttt{Trilinos}. The tests are performed on the unit ball problem for polynomial degrees $p=1,2,3$ with $s=3$ smoothing steps. Left: 2D, Right: 3D.
        For polynomial degrees higher that 1 AMG did not converge. Timings obtained on a single core of an AMD EPYC 7282. }
    \label{fig:timing_solver}

\end{figure}

Finally, we evaluate the computational efficiency of the proposed solver. Figure~\ref{fig:timing_solver} presents a
comparison of the throughput (Degrees of Freedom per second) between our multigrid preconditioner with the
Full-Residual Shy Patch smoother and an algebraic multigrid (AMG) preconditioner. Although AMG currently shows higher
throughput for $p=1$, it is worth noting that patch-based smoothers are particularly advantageous for matrix-free
implementations. Recent studies~\cite{wichrowski2025smoothers,
    cui2025multigrid,wichrowski2025local,wichrowski2025multigrid} have demonstrated that matrix-free geometric multigrid
methods with patch smoothers can achieve high performance on modern architectures by minimizing memory traffic. Our
current implementation is matrix-based, and we expect that a matrix-free optimization would yield significant
performance improvements. All timing measurements reported here were obtained using a single core of an AMD EPYC 7282
processor.

\begin{figure}[htbp]
    \centering
    \begin{tikzpicture}
    % \begin{scope}
    %     \begin{axis}[
    %             width=0.45\columnwidth,
    %             xlabel={Number of refinements}, ylabel={Ratio 2 sweeps / 1 sweep - 1}, paperplot,
    %             ymin=0,
    %             ymax=0.3 ]

    %         \immediate\write18{python3 scripts/compute_sweep_ratio.py results/full_2D.csv 3 1. > results/tmp/ratio_2D.csv}

    %         \newcommand{\addRatioPlot}[1]{% 
    %             \addplot+[thick] table[
    %                     col sep=comma,
    %                     x=refinements,
    %                     y=ratio,
    %                     restrict expr to domain={\thisrow{p}}{#1:#1},
    %                     unbounded coords=discard
    %                 ] {results/tmp/ratio_2D.csv}; \addlegendentry{$p=#1$}%
    %         }
    %         \addRatioPlot{1}
    %         \addRatioPlot{2}
    %         \addRatioPlot{3}

    %     \end{axis}
    % \end{scope}

    \begin{scope}[xshift=0.48\columnwidth]
        \begin{axis}[
                width=0.45\columnwidth,
                xlabel={Number of refinements},
                legend pos=north east,
                paperplot,
                ymin=0,
                xtick={0,1,2,3,4,5,6}]

            \immediate\write18{python3 scripts/compute_sweep_ratio.py results/full_3D.csv 4 0.0 > results/tmp/ratio_3D.csv}

            \newcommand{\addRatioPlotThree}[1]{%
                \addplot+[thick] table[
                        col sep=comma,
                        x=refinements,
                        y=ratio,
                        restrict expr to domain={\thisrow{p}}{#1:#1},
                        unbounded coords=discard ] {results/tmp/ratio_3D.csv}; \addlegendentry{$p=#1$}%
            }

            \addRatioPlotThree{1}
            \addRatioPlotThree{2}
            \addRatioPlotThree{3}

        \end{axis}
    \end{scope}
\end{tikzpicture}
    \caption{Relative increase in time per iteration when using two smoothing steps versus one, computed as $(T_{s=2}/T_{s=1})-1$, measured on the unit ball with shyness $\xi=4$ and cell threshold $\lambda=0$. Additional sweeps are restricted to boundary patches. Timings obtained on a single core of an AMD EPYC 7282.}
    %2D results are skipped here because the timings are too small to be meaningful: too much pollution and the ratio went negative for some runs, so I guess there were measurement errors. Anyway, the 3D results are more relevant.
    \label{fig:timing_sweeps}
\end{figure}
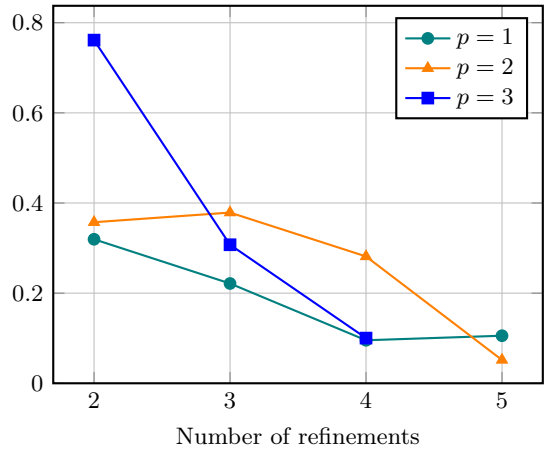

We further analyze the cost of performing additional smoothing steps. Figure~\ref{fig:timing_sweeps} displays the
relative increase in computational time when increasing the number of smoothing steps from one to two. The plotted
value corresponds to the ratio $(T_{s=2}/T_{s=1}) - 1$, representing the cost of the second (and consecutive) sweeps
relative to the first one. As detailed in Section~\ref{sec:shy_patches}, the Shy Patch smoother restricts subsequent
relaxation sweeps to only those patches affected by the boundary. Consequently, the cost of these additional sweeps is
substantially lower than that of the initial full sweep. Furthermore, as the mesh is refined, the fraction of patches
near the boundary diminishes, leading to a decreasing relative cost for the extra smoothing steps. Given the
significant reduction in iteration counts achieved with multiple smoothing steps (see Table~\ref{tab:iter_smoothing}),
this marginal additional cost is well justified, particularly for higher polynomial degrees.

\section{Conclusion}
\label{sec:conclusion}

We have addressed the algebraic challenges inherent to the Shifted Boundary Method (SBM), specifically for high-order
Continuous Galerkin discretizations. While SBM significantly simplifies mesh generation by decoupling the geometry from
the grid, it inherently shifts the complexity to the linear solver. The resulting systems exhibit condition numbers
scaling as $\mathcal{O}(h^{-2})$, which, when combined with the non-symmetry and potential indefiniteness introduced by
the boundary extrapolation, renders them resistant to standard preconditioning techniques. We have demonstrated that
this complexity can be effectively managed through a geometric multigrid preconditioner equipped with our novel
"Full-Residual Shy Patch" smoother.

The core innovation of this approach lies in the construction of local smoothing problems that are faithful
representations of the global SBM formulation. By defining subspaces over vertex patches and incorporating the full
global residual, we ensure that the non-local boundary coupling introduced by the SBM extension operators is correctly
resolved. Crucially, the introduction of the *shyness* criterion—which prevents the formation of patches around
isolated or insufficiently supported vertices—guarantees the stability of these local problems. This strategy avoids
the numerical instability associated with small cut elements that often plagues unfitted methods.

A significant advancement of this work is the method's performance at higher polynomial degrees. Unlike previous
attempts with Discontinuous Galerkin formulations where cell-wise smoothers failed to converge efficiently at $p=3$,
the Shy Patch smoother maintains consistent performance. This is particularly evident in our $p$-multigrid experiments
in 2D, where the solver demonstrated exceptional robustness with low, stable iteration counts up to polynomial degree
$p=4$. We have further shown that the computational cost can be optimized by applying additional smoothing sweeps
exclusively to boundary patches, effectively targeting the source of the error without incurring the cost of a global
sweep.

Comparisons with Algebraic Multigrid (AMG) highlight the necessity of this geometric approach; while AMG struggles
significantly with high-order SBM discretizations, our method remains stable. Although the current implementation is
matrix-based, the patch-based design is inherently compatible with tensor-product operations in the bulk domain.
Consequently, this work lays the foundation for future high-performance, matrix-free SBM solvers that can fully exploit
modern hardware architectures while retaining the geometric flexibility of unfitted methods.
\section*{Acknowledgments}
% Optional: Acknowledge funding sources, helpful discussions, etc.
The author declares the use of language models (ChatGPT, Gemini, and Claude) to improve the clarity and readability of
the manuscript. All scientific content and technical claims are solely the responsibility of the author.

\FloatBarrier
\bibliographystyle{siam} % Or any other style
\bibliography{literature} % Assuming a references.bib file will be created 
\end{document}